# REGIONAL STRATEGIC SENSORS CHARACTERIZATIONS

RAHEAM A MANSOR AL-SAPHORY [1,*], MOHAMMED A AL-JOUBORY [1] AND MAHMOOD K JASIM [2]

[1]Department of Mathematics, College of Education, Tikrit University, Iraq

[2] Department of Mathematics and Physical Sciences, College of Arts & Sciences, University of Nizwa,

Sultanate of Oman

**Abstract:** In this paper, a linear infinite dimensional distributed systems in a Hilbert space has been discussed and analyzed where the dynamics of system is governed by strongly continuous semi-groups. The characterizations of regional strategic sensors have been given and tackled for different cases of regional observability. Furthermore, the results so obtained are applied to two-dimensional systems. Various cases of sensors are considered and analyzed. Also, the authors show that, the existent of a sensor for the diffusion system is not strategic in the usual sense, but it may be regionally strategic of this system.

## 1. Introduction

Many works, in distributed parameter systems (DPSs), have been devoted to the observation problem [1]. It has often been studied independently of any geometric considerations, and most of the works were focused on the observation and reconstruction of the state in a certain observation space [2-3]. The notion of sensors and actuators introduced in the 1980s by El Jai and Pritchard allows for a better description of measurements and actions [4]. In addition, the study of observability and controllability can be considered with respect to the structure, number and location of sensors and actuators [5-6]. For linear DPSs, observability and controllability are dual notions and most results on observability can be deduced from those on controllability by duality [7-8]. The regional observability concept has been developed recently by El Jai *et al.*, in the 1993 [8]. In this case of regional analysis, the problem concerns the state observation in a sub-region $\omega$ of system domain $\Omega$ [10-11].

---

[*]Corresponding author

The study of this notion motivated by certain concrete-real problem, in thermic, mechanic, environment [11-14]. If a system is defined on a domain $\Omega$ and represented by the room model as in the figure 1below, then we are interested in the regional state on $\omega$ of the domain $\Omega$ [15].

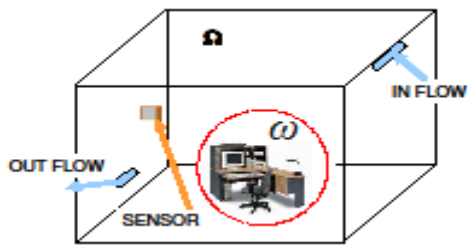

*Fig. (1) Room model of sensor and actuator problem.*

The purpose of this paper is to give sufficient conditions of strategic sensor in this region which observes regional state. This paper is organized as follows:

The second section is focused on the considered system and the problem of regional observability. The third section is devoted to the mathematical concept of regional observability and the characterization of regional strategic sensors in various situations is studied. In the last section, we illustrate applications with many situations of sensor locations.

## 2. Regional observability

In this section, we are interested to study the notion of regional observability and recall original results related to particular systems.

**2.1 Problem statement**

Let $\Omega$ be a regular bounded open subset of $R^n$, with boundary $\partial\Omega$ and $[0,T]$, $T>0$ be a time measurement interval. Suppose that $\omega$ be a non-empty given sub-region of $\Omega$. We denote $\Theta = \Omega \times ]0,T[$ and $\Pi = \partial\Omega \times ]0,T[$. The considered distributed parabolic systems is described by the following state space equations

$$\begin{cases} \dfrac{\partial x}{\partial t}(\xi,t) = Ax(\xi,t) + Bu(t) & \Theta \\ x(\xi,0) = x_\circ(\xi) & \Omega \\ x(\eta,t) = 0 & \Pi \end{cases} \quad (1)$$

augmented with the output function

$$y(.,t) = C x(.,t) \quad (2)$$

where $A$ is a second order linear differential operator, which generates a strongly continuous semi-group $(S_A(t))_{t\geq 0}$ on the Hilbert space $X = L^2(\Omega)$ and is self-adjoin with compact resolvant. The operator $B \in L(R^p, X)$ and $C \in L(R^q, X)$, depend on the structures of actuators and sensors [16-18]. The spaces $X, U$ and $O$ be separable Hilbert spaces where $X$ is the state space, $U = L^2(0,T,R^p)$ is the control space and $O = L^2(0,T,R^q)$ is the observation space, where $p$ and $q$ are the numbers of actuators and sensors. Under the given assumption, the system (1) has a unique solution [1]:

$$x(\xi,t) = S_A(t)x_\circ(\xi) + \int_0^t S_A(t-\tau)Bu(\tau)\,d\tau \quad (3)$$

The problem is that, how to present sufficient conditions for regional strategic sensors which enable to observe the current state in a given sub region $\omega$ (see figure 1above), using convenient sensors.

### 2.2 Definitions and characterizations

The regional observability concept has been developed recently by El Jai [9-10] and extended to the regional asymptotic state by Al-Saphory and El Jai in ref.s [13,16-18]. To recall regional observability, consider the associated autonomous system to (1) given by

$$\begin{cases} \dfrac{\partial x}{\partial t}(\xi,t) = Ax(\xi,t) & \Theta \\ x(\xi,0) = x_\circ(\xi) & \Omega \\ x(\eta,t) = 0 & \Pi \end{cases} \quad (4)$$

With $x(\xi,0)$ is supposed known in $L^2(\Omega)$. Thus, the knowledge of $x(\xi,0)$ permit to observe regional state $x(\xi,t)$ at any time $t$.

- The measurements are obtained from the output function (2). In this case, the

solution of (4) is given by the following form,

$$x(\xi,t) = S_A(t)x_\circ(\xi), \quad \forall t \in [0,T] \tag{5}$$

- We define the operator

$$K : x \in X \to Kx = CS_A(.)x \in O \tag{6}$$

Then, we obtain

$$y(.,t) = K(t)x(.,0)$$

where $K$ is bounded linear operator (this is valuable on some output function) [19].

- We note that $K^* : O \to X$ is the adjoint operator of $K$ defined by

$$K^* y^* = \int_0^t S_A^* C^* y^* ds \tag{7}$$

- Consider a sub-region $\omega \subset \Omega$ and let $\chi_\omega$ is the restriction function defined by

$$\chi_\omega : X = L^2(\Omega) \to L^2(\omega) \tag{8}$$

$$x \to \chi_\omega x = x|_\omega$$

where $x|_\omega$ is the restriction of $x$ to $\omega$.

- The adjoint $\chi_\omega^* : L^2(\omega) \to L^2(\Omega)$ of $\chi_\omega$ is defined by

$$(\chi_\omega^* x)\xi = \begin{cases} x(\xi) & \xi \in \omega \\ 0 & \xi \in \Omega \setminus \omega \end{cases}$$

Now, to present strategic sensors, we need some definitions and characterizations of regional observability concept as in [11-12].

**Definition 2.1.** The system (4) augmented with the output function (2) is said to be regionally exactly observable on $\omega$ (or exactly $\omega$- observable), if

$$\operatorname{Im} \chi_\omega K^* = L^2(\omega)$$

**Definition 2.2.** The system (4) augmented with the output function (2) is said to be regionally approximately observable on $\omega$ (or approximately $\omega$- observable), if

$$\overline{\operatorname{Im} \chi_\omega K^*(.)} = L^2(\omega)$$

**Remark 2.3.** *The definition 2.2. is equivalent to say that the system (4)-(2) is approximately $\omega$-observable if*

$$\ker \quad K(t)\chi_\omega^* = \{0\}$$

Then, the following characterizations can extend to the regional case as in ref. [7].

**Proposition 2.4.** *The system (4)-(2) is exactly $\omega$-observable if there exists $\nu > 0$ such that $\forall x_\circ \in L^2(\omega)$,*

$$\|\chi_\omega x_\circ\|_{L^2(\omega)} \leq \nu \|K\chi^*_\omega x_\circ\|_{L^2(0,T,O)} \tag{9}$$

**Proof.** The proof of this property is deduced from the usual results on observability considering $\chi_\omega K^*$ [3].

Let $E, F$ and $G$ be Banach reflexive space and $f \in L(E, G), g \in L(F, G)$, then we have

(1) $Im\ f \subset Im\ g$
(2) then there exist $c > 0$ such that

$$\|f^* x^*\|_{E^*} \leq \|g^* x^*\|_{F^*}, \forall\ x^* \in G^*.$$

Now, if this result is applied. Choosing

$$E = G = L^2(\omega), \quad F = 0, f = Id_{L^2(\omega)}$$

and

$$g = \chi_\omega K^*,$$

therefore, we obtain the inequality (9) ∎.
From the proposition 2.4 .we can get the following result.

**Corollary 2.5.** *we have:*

*(1) The notion of approximate $\omega$-observability is far less restrictive than the exact $\omega$-observability.*

*(2) From the equation (9) there exists a reconstruction error operator that gives an estimation $\tilde{x}_\circ$ of the initial state $x_\circ$ in $\omega$ [10]. Then, we have*

$$\left\| x_\circ - \tilde{x}_\circ \right\|_{L^2(\omega)} \leq \left\| x_\circ - \tilde{x}_\circ \right\|_{L^2(\Omega)} \tag{10}$$

**Remark 2.6.** *The regional observability concept is more convenient in the analysis of real systems. We can deduce that:*

*(1) The definitions 2.1. and 2.2. are general and can be applied to the case where $\omega = \Omega$.*

*(2) The equation (10) shows that he regional state reconstruction will be more precise than if we estimate the state in the whole domain $\Omega$.*

*(3) A system is exactly observable, then, it is exactly $\omega$-observable, but the converse is not true in general. Now, we prove that property (3) of remark 2.6.*

**Proof.** We see that if the system is exactly observable on $\Omega$, then it is exactly $\omega$-observable and this is a consequence of (10) and then

$$\| x_\circ \|_{L^2(\omega)} \leq \| x_\circ \|_{L^2(\Omega)} \qquad \forall x_\circ \in L^2(\omega)$$

We can explain this by:

$$\|x_\circ\|_{L^2(\omega)} = \int_\omega |x_\circ|^2 < \infty$$

and

$$\|x_\circ\|_{L^2(\Omega)} = \int_\Omega |x_\circ|^2 < \infty$$

since $L^2(\Omega) = \left\{x_\circ : \int_\Omega |x_\circ|^2 < \infty\right\}$ and $L^2(\omega) = \left\{x_\circ : \int_\omega |x_\circ|^2 < \infty\right\}$

Here, if $\omega \subset \Omega$, then $|x_\circ|_\omega \leq |x_\circ|_\Omega$

$$\Rightarrow \int_\omega |x_\circ|_\omega^2 \leq \int_\Omega |x_\circ|_\Omega^2$$

$$\Rightarrow \|x_\circ\|_{L^2(\omega)} \leq \|x_\circ\|_{L^2(\Omega)}$$

$$\therefore \Rightarrow \|\chi_\omega x_\circ\|_{L^2(\omega)} \leq \|x_\circ\|_{L^2(\Omega)} \leq \|Kx_\circ\|_{L^2(0,T,O)}$$

Then from proposition 2.6., then system (4)-(2) is exactly $\omega$-observable.

## 3. $\omega$-Strategic sensors

The purpose of this section is to give the characterization for sensors in order that the system (1) is regionally observable in a region $\omega$.

### 3.1 Sensors notion

This subsection consists of the concept of the sensors, which was coined by A. El Jai [4]. Sensors can play a fundamental role in the understanding of any real system. Sensors form an important link between a system and its environment. Sensors have a passive role and allow the system evolution to be measured. This means that the sensors will define an output function. In any case of sensors is considered via a space variable, mathematically speaking, the space variable is present in all systems described by partial differential equations. Now, let $\Omega$ be the spatial domain in which the system (4) is modeled. We can have much geometry for the support given by the following cases (see figure 2 below):

(1) Pointwise or zone sensors inside the physical domain $\Omega$.
(2) Pointwise or zone sensors on the boundary of the physical domain $\Omega$.
(3) Filament sensors in the domain $\Omega$.
(4) Mobile pointwise sensors in the domain $\Omega$.

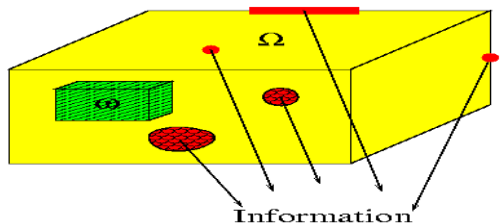

*Fig. (2) The domain of $\Omega$, the sub-region $\omega$, various sensors locations.*

**Definition 3.1.** A sensor may be defined by any couple $(D, f)$ where:

(1) $D$ denotes a closed subset of $\overline{\Omega}$, which is the support of sensor in the zone case $D = \{b\}$ in the pointwise case),

(2) $f \in L^2(\Omega)$ defines the spatial distribution of the sensor in the zone case ( $f = \delta_b$ in the pointwise case, $\delta_b$ is a Dirac mass in $b$ ).

**Remark 3.2.** *According to the choice of the parameters $D$ and $f$ we have various types of sensors. In this case of $q$ sensors, we consider $(D_i, f_i)_{1 \leq i \leq q}$ with*

$$D_i \subset \Omega \ (or \ D_i \subset \partial\Omega \ )$$

*and*

$$f_i \in L^2(D_i), \ D_i \cap D_j = \phi, \ if \ 1 \leq i \neq j \leq q.$$

*The associated output function is:*

$$y(t) = \left[ y_1(t), y_2(t), ..., y_q(t) \right]^T$$

*where*

$$y_i(t) = < x(t), f_i >_{L^2(D_i)}$$

*In the case of zone sensors and then,*

$$y(t) = [< x(t), f_1 >_{L^2(D_1)}, ..., < x(t), f_q >_{L^2(D_q)}]^T \tag{11}$$

*From these equations the output function may be written in the form*

$$y(t) = C\, x(t)$$

*If the sensors (or some of them) are pointwise, the correspondent components to the associated output are given by:*

$$y(t) = \left[ y_1(t), y_2(t), ..., y_q(t) \right]^T \tag{12}$$

$$= < x(t), \delta_{b_i} >_{L^2(\Omega)}$$

$$= x(b_i, t)$$

*where $\delta_{b_i}$ denotes the Dirac mass concentrated in $b_i$. Now, let us consider the orthonormal set of eigenfunctions ($\varphi_{nj}$) of A is associated with the eigenvalues ($\lambda_n$)* of multiplicity $r_n$, then, the state vector $x(t)$ can be expressed by the following form:

$$x(t) = \sum_{n,j=1}^{\infty} \varphi_{nj} x_n(t) \tag{13}$$

and we have the following forms:

(i) **Zone sensors case:** In this case, the operator $C$ is of type $q \times \infty$ with

$$C_{ij} = <\varphi_{nj}, f_i>_{L^2(D_i)}$$

$$= \int_{D_i} \varphi_{nj}(\xi) f_i(\xi) d\xi, \ 1 \leq i \leq q, \ 1 \leq j \leq r_n \ \text{and} \ n \geq 1$$

and then the output function can be given by

$$y(t) = \begin{bmatrix} \sum_{n \geq 1} x_n(t) \sum_{j=1}^{r_n} <\varphi_{nj}, f_1>_{L^2(D_1)} \\ . \\ . \\ . \\ \sum_{n \geq 1} x_n(t) \sum_{j=1}^{r_n} <\varphi_{nj}, f_q>_{L^2(D_q)} \end{bmatrix} \tag{14}$$

Then from (13) the output function (12) can be written in the form

$$y(.,t) = Cx(.,t) = \int_{D_i} x(\xi,t) f_i(\xi) d\xi \tag{15}$$

**Pointwise sensors case:** The operator $C$ is of type $q \times \infty$ with

$$C_{ij} = <\varphi_{nj}, \delta_{b_i}>_{L^2(\Omega)} = \varphi_{nj}(b_i)$$

$$= \int_{\Omega} \varphi_{nj}(\xi,t) \delta_{b_i}(\xi - b_i) d\xi, \ 1 \leq i \leq q, \ 1 \leq j \leq r_n \ \text{and} \ n \geq 1$$

Thus we have

$$y(t) = \begin{bmatrix} \sum_{n \geq 1} x_n(t) \sum_{j=1}^{r_n} <\varphi_{nj}, \delta_{b_i}>_{L^2(\Omega)} \\ . \\ . \\ . \\ \sum_{n \geq 1} x_n(t) \sum_{j=1}^{r_n} <\varphi_{nj}, \delta_{bq}>_{L^2(\Omega)} \end{bmatrix} \tag{16}$$

In this case $C$ is not bounded operator in $X$ [5]. Then, from above the output function (12) can be given by the form

$$y(.,t) = Cx(.,t) = \int_\Omega x(\xi,t)\delta_{b_i}(\xi - b_i)\,d\xi \tag{17}$$

**Boundary sensors case:** In the case of boundary measurements (pointwise or zone) the support of sensors $D_i$ is subset of $\partial\Omega$.

i.e. the output function can be obtained by

$$y(t) = \begin{bmatrix} \sum_{n \geq 1} x_n(t) \sum_{j=1}^{r_n} < \dfrac{\partial \varphi_{nj}}{\partial \nu}, \delta_{b_1} >_{L^2(\partial\Omega)} \\ . \\ . \\ . \\ \sum_{n,j} x_n(t) \sum_{j=1}^{r_n} < \dfrac{\partial \varphi_{nj}}{\partial \nu}, \delta_{b_q} >_{L^2(\partial\Omega)} \end{bmatrix} \tag{18}$$

and

$$y(.,t) = Cx(.,t) = \int_\Omega \dfrac{\partial x}{\partial \nu}(\xi,t)\,\delta_{b_i}(\xi - b_i)\,d\xi \tag{19}$$

Now in the case where the zone measurements, with $D_i = \Gamma_i \subset \partial\Omega$ and $f_i \in L^2(\Gamma_i)$, then the output function can be written by

$$C_{ij} = < \dfrac{\partial \varphi_{nj}}{\partial \nu}, f_i >_{L^2(\Gamma_i)}$$

$$= \int_{\Gamma_i} \dfrac{\partial \varphi_{nj}}{\partial \nu}(\xi) f_i(\xi)\,d\xi, \quad 1 \leq i \leq q,\; 1 \leq j \leq r_n \text{ and } n \geq 1$$

and then

$$y(t) = \begin{bmatrix} \sum_{n \geq 1} x_n(t) \sum_{j=1}^{r_n} < \dfrac{\partial \varphi_{nj}}{\partial \nu}, f_1 >_{L^2(\Gamma_1)} \\ . \\ . \\ . \\ \sum_{n \geq 1} x_n(t) \sum_{j=1}^{r_n} < \dfrac{\partial \varphi_{nj}}{\partial \nu}, f_q >_{L^2(\Gamma_q)} \end{bmatrix} \tag{20}$$

The operator $C$ is not bounded in this case [5]. Then, the output function (12) can be given by the form

$$y(.,t) = Cx(.,t) = \int_{\Gamma_i} \dfrac{\partial x}{\partial \nu}(\xi,t)\,f_i(\xi)\,d\xi \tag{21}$$

**Definition 3.3.** A sensor $(D, f)$ is $\omega$-strategic if the corresponding system (14) augmented with the output function (2) is approximately $\omega$-observable.

**Definition 3.4.** A suite of $(D_i, f_i)_{1 \leq i \leq q}$ is said to be $\omega$-strategic if there exist at least one sensor $(D_1, f_1)$ which is approximately $\omega$-strategic.

We can deduce that the following result:

**Corollary 3.5.** *A sensor is $\omega$-strategic if the corresponding system (12)-(14) is exactly $\omega$-observable.*

**Proof.** Let the system (12)-(14) is exactly $\omega$-observable. Then, we have

$$\operatorname{Im} \chi_\omega K^* = L^2(\omega)$$

From the decomposition sub-spaces of direct sum in Hebert space, we can represent $L^2(\Omega)$ by the unique form [3]

$$\ker \chi_\omega + \operatorname{Im} K^* = L^2(\Omega)$$

we obtain

$$K(t)\chi_\omega^* = \{0\}$$

This is equivalent to

$$\overline{\operatorname{Im} \chi_\omega K^*(.)} = L^2(\omega)$$

Finally, we can deduce this system is approximately $\omega$-observable and therefore this sensor is $\omega$-strategic. ∎

Thus, the definition 2.1., proposition 2.4. and corollary 3.5. guarantee $\omega$-strategic sensors with far more restrictive conditions.

**Remark 3.6.** *From the previous results, we note that:*
*(1) a sensor which is strategic for a system, is $\omega$-strategic.*
*(2) a sensor which is $\omega_1$-strategic for a system where $\omega_1 \subset \Omega$, is $\omega_2$-strategic for any $\omega_2 \subset \omega_1$.*
*(3) One can find various sensors which are not strategic in usual sense for systems, but may be $\omega$-strategic and achieve the observability in $\omega$. This is illustrated in the following counter-example.*

### 3.2 A counter-example
Consider the diffusion process

$$\begin{cases} \dfrac{\partial x}{\partial t}(\xi,t) = \dfrac{\partial^2 x}{\partial \xi^2}(\xi,t) & \Omega \times ]0,1[ \\ x(0,t) = x(1,t) = 0 & \partial\Omega \\ x(\xi,0) = x_\circ(\xi) & \Omega \end{cases} \quad (22)$$

where $\Omega = ]0,1[$ and the function output

$$y(.,t) = \int_\Omega x(\xi,t)\delta(\xi - b)\,d\xi = x(b,t), \quad b \in \Omega \quad (23)$$

where $b$ is the location of pointwise sensor $(\delta_b, b)$ in $\Omega$ as defined in figure 3

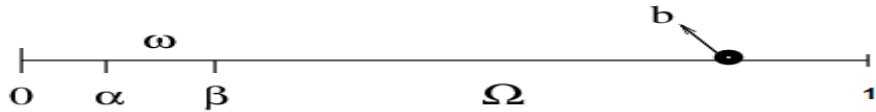

*Fig. (3) The domain $\Omega$, the sub-region $\omega$ and the location sensors $b$.*

First, we must prove that the system (22)-(23) is not approximately observable in $\Omega$, that means the sensor $(\delta_b, b)$ is not strategic. For this purpose, we can write the system (22) as a state space one dimensional system

$$\dot{x}(\xi,t) = Ax(\xi,t) \qquad (24)$$
$$x(0) = x_\circ(\xi)$$

Where $A = \dfrac{\partial^2}{\partial \xi^2}$ generates the continuous semi-group $(S(t))_{t \geq 0}$ given by [7].

$$S(t)x_\circ = \sum_{i=1}^{\infty} \exp(\lambda_i t) <x_\circ, \varphi_i>_{L^2(\Omega)} \varphi_i$$

Where, $\varphi_i, \lambda_i$ are the eigenfunctions and the associated eigenvalues of $A$. From the solution of (22), we have

$$y(\xi,t) = \sum_{i=1}^{\infty} \exp(\lambda_i t) <x_\circ, \varphi_i>_{L^2(\Omega)} \varphi_i(b) = CS(t)x_\circ = K(t)x_\circ$$

then, the system (22)-(23) is approximately observable if

$$\text{Ker } K(t) = \{0\}$$

As proved by El Jai and Pritchard [7], if $b \in Q$ (the rational number), then the system (22)-(23) is not approximately observable and the sensor $(\delta_b, b)$ is not strategic. This is the case if we consider

$$J = \{j \mid jb \in 2N\}$$

then

$$\text{Ker } K(t) = \overline{\{\varphi_i\}}_{j \in J}$$

Now, we can prove that for certain $j_\circ \in J$, the state $\varphi_{j_\circ}$ that is not approximately observable in $\Omega$ but can be approximately $\omega$-observable on certain sub-region

$\omega \subset \Omega$. Let

$\omega = [\alpha, \beta] \subset [0,1]$, with $\beta = \alpha + b$, we have

$$\int_\alpha^\beta \varphi_j^2 d\xi = \beta - \alpha \quad \forall j \in J \quad \text{and for} \quad i, j \in J, i \neq j,$$

we obtain

$$\int_\alpha^\beta \varphi_i \varphi_j d\xi = 0$$

If $j_\circ \in J$ then $\varphi_{j_\circ}$ is not approximately observable on $[0,1]$ and the sensor is not strategic. Let us show that $\varphi_{j_\circ}$ is approximately $\omega$-observable on $\omega = [\alpha, \beta]$. we have

$K(t) \chi_\omega^* \chi_\omega \varphi_\circ = 0, \forall t \in [0, T]$. Then

$$\sum_{i \notin j} \exp(\lambda_i t) < \varphi_{j_\circ}, \varphi_i >_{L^2(\omega)} \varphi_i(b) = 0$$

$$\Rightarrow < \varphi_{j_\circ}, \varphi_i >_{L^2(\omega)} = 0 \quad \forall i \notin J$$

for $i = i_\circ$ such that $i_\circ b = 2k + 1, k \in N$, we have

$$\Rightarrow < \varphi_{j_\circ}, \varphi_i >_{L^2(\omega)} = 0 \Leftrightarrow i_\circ \tan j_\circ \pi\alpha = j_\circ \tan i_\circ \pi\alpha, \quad \forall i_\circ, j_\circ$$

which is not true, in general. Consider the space where $\alpha = \frac{1}{4}$, $b = \frac{1}{2}$, $i_\circ = 6$, $j_\circ = 4$. Then $x_\circ = \sqrt{2} \sin 6\pi y$ is approximately $\omega$-observable on $\left[\frac{1}{4}, \frac{3}{4}\right]$ and the sensor is $\omega$-strategic. ■

The concept of regional strategic sensor on $\omega$ can be characterized by the following results:

**Theorem 3.7.** *Assume that* $\sup r_n = r < \infty$*, then the suite of sensors* $(D_i, f_i)_{1 \leq i \leq q}$ *is* $\omega$- *strategic if and only if*

  *(1)* $q \geq r$

  *(2) rank* $G_n = r_n$

*where*

$$G_n = (G_n)_{ij} \qquad 1 \leq i \leq q, \quad 1 \leq j \leq r_n$$

*and*

$$(G_n)_{ij} = \begin{bmatrix} <\varphi_{n_1}, f_1(\,.\,)>_{L^2(\Omega_1)}, \ldots, <\varphi_{n_{r_n}}, f_1(\,.\,)>_{L^2(\Omega_1)} \\ \vdots \\ <\varphi_{n_1}, f_q(\,.\,)>_{L^2(\Omega_p)}, \ldots, <\varphi_{n_{r_n}}, f_q(\,.\,)>_{L^2(\Omega_p)} \end{bmatrix}$$

**Proof.** The proof is developed in the case where the sensors are of zone type and located inside the domain $\Omega$. If suite of sensors are $\omega$-strategic, then the corresponding system (22)-(24) is approximately $\omega$-observable, it is equivalent to

$$[K\chi_\omega^* x^* = 0 \Rightarrow z^* = 0], \text{ for } x^* \in L^2(\omega), [6] \text{ we have}$$

$$K\chi_\omega^* x^* = \left( \sum_{n \geq 1} \exp(\lambda_n t) \sum_{j=1}^{r_n} <\varphi_{nj}, x^*>_\omega <\varphi_{nj}, f_i> \right), \quad 1 \leq i \leq q$$

If the suite of sensors is not strategic sensors, i.e, the system (22)-(24) is not approximately $\omega$-observable, then there exists $x^* \neq 0$ such that

$$K\chi_\omega^* x^* = 0 \Leftrightarrow \sum_{j=1}^{r_n} <\varphi_{nj}, x^*>_\omega <\varphi_{nj}, f_i> = 0 \quad \forall n, \ n \geq 1$$

Let $x_n$ define by

$$x_n = \begin{bmatrix} <\varphi_{n_1}, x^*>_{L^2(\omega)} \\ \ldots \\ <\varphi_{n\,r_n}, x^*>_{L^2(\omega)} \end{bmatrix}$$

Then,

$$G_n x_n = 0, \quad \forall n \geq 1 \Leftrightarrow rank\, G_n \neq r_n.$$

Now, we to prove conversely if $rank\ G_n \neq r_n$ for some $n$, then there exists

$$x_n = \begin{bmatrix} x_{n_1} \\ \vdots \\ x_{n_{r_n}} \end{bmatrix} \neq 0, \qquad x^* = \sum_{j=1}^{r_n} x_{nj} \varphi_{nj} \in L^2(\omega) \neq 0$$

such that

$$G_n x_n = 0$$

so we can construct a non zero $x^* \in L^2(\omega)$. Considering

$$<x^*, \varphi_{jk}>_{L^2(\omega)} = 0 \quad \text{if} \quad j = n \quad \text{and}$$

$$<x^*, \varphi_{nk}>_\omega = x_{nk}, \quad 1 \leq k \leq r_n$$

for which
$$\sum_{k=1}^{r_j} <\varphi_{jk}, f_i>_{\Omega_i} <x^*, \varphi_{jk}> = 0 \quad j \neq n, \ 1 \leq i \leq q \quad \text{and also}$$

$$\sum_{k=1}^{r_j} <\varphi_{nk}, f_i>_{\Omega_i} <x^*, \varphi_{nk}> = 0 \quad 1 \leq i \leq q$$

otherwise there exists $x^* \neq 0 \in L^2(\omega)$, such that

$$K\chi_\omega^* x^* = 0,$$

Thus, the system (22)-(24) is not approximately $\omega$-observable and then the sensors are not $\omega$-strategic. ∎

**Corollary 3.8.** *If the system (22)-(24) is exactly $\omega$-observable, rank condition in theorem 3.7. is satisfied.*

**Remark 3.9.** *The previous result can be extended to the case of pointwise, filament sensors as in ref.s [16-18].*

## 4. Application to sensor locations

In this section, we present an application of the above results to a two-dimensional system defined on $\Omega = ]0, 1[ \times ]0, 1[$ by the form

$$\begin{cases} \dfrac{\partial x}{\partial t}(\xi_1, \xi_2, t) = \Delta x(\xi_1, \xi_2 t) & \Theta \\ x(\xi_1, \xi_2, 0) = x_\circ(\xi_1, \xi_2) & \Omega \\ x(\eta_1, \eta_2, t) = 0 & \Pi \end{cases} \quad (25)$$

together with output function by (2). Let $\omega = ]\alpha_1, \beta_1[ \times ]\alpha_2, \beta_2[$ be the considered region is subset of $]0, 1[ \times ]0, 1[$. In this case, the eigenfunctions of system (25) are given by

$$\varphi_{ij}(\xi_1, \xi_2) = \dfrac{2}{\sqrt{(\beta_1 - \alpha_1)(\beta_2 - \alpha_2)}} \sin i\pi(\dfrac{\xi_1 - \alpha_1}{\beta_1 - \alpha_1}) \sin j\pi(\dfrac{\xi_2 - \alpha_2}{\beta_2 - \alpha_2}) \quad (26)$$

associated with eigenvalues

$$\lambda_{ij} = -\left(\dfrac{i^2}{(\beta_1 - \alpha_1)^2} + \dfrac{j^2}{(\beta_2 - \alpha_2)^2}\right) \quad (27)$$

The following results give information on the location of internal zone or pointwise $\omega$-strategic sensors.

### 4.1 Internal zone sensor

Consider the system (25) together with output function (2) where the sensor supports $D$ are located in $\Omega$. The output (2) can be written by the form

$$y(t) = \int_D x(\xi_1, \xi_2, t) f(\xi_1, \xi_2) d\xi_1 d\xi_2 \tag{28}$$

where $D \subset \Omega$ is location of zone sensor and $f \in L^2(D)$. In this case of (figure 4), the eigenfunctions and the eigenvalues

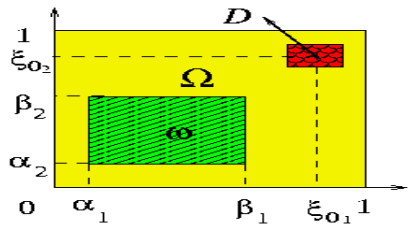

*Fig. (4) Domain $\Omega$, sub-region $\omega$ and location $D$ of internal zone sensor.*

are given by (4.2) and (4.3). However, if we suppose that

$$\frac{(\beta_1 - \alpha_1)^2}{(\beta_2 - \alpha_2)^2} \notin \mathbb{Q} \tag{29}$$

Then $r = 1$ and one sensor $(D, f)$ my be sufficient to achieve $\omega$-observability of systems (25)-(28) [20]. Let the measurement support is rectangular with

$$D = [\xi_1 - l_1, \xi_1 + l_2] \times [\xi_2 - l_2, \xi_2 + l_2] \in \Omega$$

then, we have the following result:

**Corollary 4.1.** *If $f_1$ is symmetric about $\xi_1 = \xi \circ_1$ and $f_2$ is symmetric about $\xi_2 = \xi_{\circ_2}$, then the sensor $(D, f)$ is $\omega$-strategic if*

$$\frac{i(\xi_{\circ_1} - \alpha_1)}{(\beta_1 - \alpha_1)} \quad \text{and} \quad \frac{i(\xi_{\circ_2} - \alpha_2)}{(\beta_2 - \alpha_2)} \notin \mathbb{N} \quad \text{for some } i.$$

### 4.2 Internal pointwise sensor

Let us consider the case of pointwise sensor located inside of $\Omega$. The system (25) is augmented with the following output function:

$$y(t) = \int x(\xi_1, \xi_2, t) \delta(\xi_1 - b_1, \xi_2 - b_2) d\xi_1 d\xi_2 \tag{30}$$

where $b = (b_1, b_2)$ is the location of pointwise sensor as defined in (figure 5)

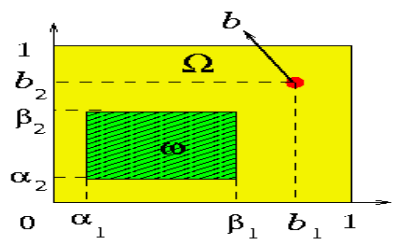

**Fig. (5) Rectangular domain, and location $b$ of internal pointwise sensor.**

If $(\beta_1 - \alpha_1)/(\beta_2 - \alpha_2) \notin Q$ then $m = 1$ and one sensor $(b, \delta_b)$ may be sufficient for $\omega$-observability of systems (25)-(30).

Thus, we obtain the following result,

**Corollary 4.2.** *The sensor $(b, \delta_b)$ is $\omega$-strategic if $i(b_1 - \alpha_1)/(\beta_1 - \alpha_1)$ and $i(b_2 - \alpha_2)/(\beta_2 - \alpha_2) \in N$, for every $i$.*

### 4.3 Internal filament sensor

Consider the case where the information is given on the curve $\sigma = \text{Im}(\gamma)$ with $\gamma \in C^1(0,1)$ (see figure 6), then we have,

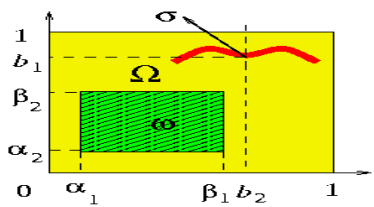

**Fig. (6) Rectangular domain, and location $\sigma$ of internal filament sensors.**

**Corollary 4.3.** *If the measurements recovered by filament sensor $(\sigma, \delta_\sigma)$ such that is symmetric with respect to the line $\xi = \xi_\circ$. Then the sensor $(\sigma, \delta_\sigma)$ is $\omega$-strategic if $i(\xi_{\circ 1} - \alpha_1)/(\beta_1 - \alpha_1)$ and $i(\xi_{\circ 2} - \alpha_2)/(\beta_2 - \alpha_2) \in N$ for all $i = 1,...,q$.*

**Remark 4.4.** *These results can be extended to the following:*

    *1. Case of Neumann or mixed boundary conditions [4-5].*

    *2. Case of disc domain $\Omega = (D,1)$ and $\omega = (0, r_\omega)$ where $\omega \subset \Omega$ and $0 < r_\omega < 1$ [16-17].*

    *3. Case of boundary sensors where $C \notin L(X, R^q)$, we refer to see [19, 21].*

    *4. we can show that the observation error decreases when the number and support of sensors increases [22].*

## 5. Conclusions

The regional strategic sensors have been developed. A various regional observability have been discussed and analyzed which permit us to avoid some bad sensor locations. Various interesting results concerning the choice of such sensors are given and illustrated in specific situations with diffusion systems. Many questions still opened, for example, the simulations of this model are under consideration and the problem of finding an optimal sensor location ensuring such an objective.